\documentclass[12pt]{amsart}
\usepackage{amssymb}
\usepackage{amsmath}
\numberwithin{equation}{section}
\usepackage{color}
\usepackage{graphicx}
\usepackage{epstopdf}

\newtheorem{thm}{Theorem}[section]
\newtheorem{lemma}[thm]{Lemma}
\newtheorem{cor}[thm]{Corollary}
\newtheorem{prop}[thm]{Proposition}
\newtheorem{con}[thm]{Conjecture}

\theoremstyle{definition}

\newtheorem{prob}[thm]{Problem}

\newtheorem{eg}[thm]{Example}

\theoremstyle{remark}

\renewcommand{\S}{\mathfrak S}
\newcommand{\s}{\sigma}

\newcommand\C{{\mathbb{C}}}

\newcommand\Q{{\mathbb Q}}
\newcommand\Z{{\mathbb{Z}}}

\newcommand\PP{{\mathbb{P}}}
\newcommand\N{{\mathbb{N}}}

\newcommand\bq{\begin{equation}}
\newcommand\eq{\end{equation}}
\newcommand\beq{\begin{eqnarray*}}
\newcommand\eeq{\end{eqnarray*}}
\newcommand\ben{\begin{enumerate}}
\newcommand\een{\end{enumerate}}
\newcommand\bit{\begin{itemize}}
\newcommand\eit{\end{itemize}}

\newcommand\des{{\rm des}}
\newcommand\exc{{\rm exc}}
\newcommand\inv{{\rm inv}}
\newcommand\maj{{\rm maj}}
\newcommand\Inv{{\rm INV}}

\newcommand\sg{{\mathfrak S}}
\newcommand\Des{{\rm DES}}
\newcommand\Dex{{\rm DEX}}
\newcommand\Exc{{\rm EXC}}

\newcommand\ch{{\rm ch}}

\newcommand\x{{\mathbf x}}

\newcommand\Par{{\rm Par}}
\newcommand\pa{{\rm par}}

\newcommand\rmaj{{\rm rmaj}}
\newcommand\asc{{\rm asc}}

\newcommand\sym{\mathfrak S}

\def\ov{\overline}

\def\zz{{\mathbb Z}}
\def\nn{{\mathbb N}}

\def\cc{{\mathbb C}}
\def\pp{{\mathbb P}}

\def\cq{{\mathcal Q}}
\def\ci{{\mathcal I}}

\def\inc{{\rm{inc}}}

\def\bft{s}
\def\hess{{\mathcal H}}

\def\Ps{{\rm \bf ps}}

\def\v{{\mathcal V}}

\begin{document}

\title[Chromatic quasisymmetric functions]
{ Chromatic quasisymmetric functions and Hessenberg varieties}
\author[Shareshian]{John Shareshian$^1$}
\address{Department of Mathematics, Washington University, St. Louis, MO 63130}
\thanks{$^{1}$Supported in part by NSF Grant
 DMS 0902142}
\email{shareshi@math.wustl.edu}

\author[Wachs]{Michelle L. Wachs$^2$}
\address{Department of Mathematics, University of Miami, Coral Gables, FL 33124}
\email{wachs@math.miami.edu}
\thanks{$^{2}$Supported in part by NSF Grant
DMS 0902323}

\begin{abstract} We discuss three distinct topics of independent interest; one in enumerative combinatorics, one in symmetric function theory, and one in algebraic geometry.
   The topic in enumerative combinatorics concerns a q-analog of a generalization of the Eulerian polynomials, the one in symmetric function theory deals with a refinement of the chromatic symmetric functions of Stanley, and the one in algebraic geometry deals with Tymoczko's representation of the symmetric group   on the cohomology of the regular semisimple Hessenberg variety of type A.
Our purpose is to explore some  remarkable connections between these topics.\end{abstract}


\date{June. 11, 2011;  version July 6, 2012}

\maketitle

\tableofcontents

\section{Introduction}

Let $H(z) :=\sum_{n \ge 0} h_n z^n$, where  $h_n$ is the complete homogeneous symmetric function of degree $n$.  The formal power series   \bq \label{hsymeq} {(1-t) H(z) \over H(tz) - t H(z)},\eq
has occurred in various contexts in the literature, including
\begin{itemize}
\item the study of a class of q-Eulerian polynomials
\item the enumeration of Smirnov words by descent number
\item the study of the toric variety associated with the Coxeter complex of the symmetric group.
\end{itemize}

We briefly describe these contexts now  and give  more detail in Section~\ref{eultorsec} with definitions in Section~\ref{prelimsec}.   In \cite{ShWa1,ShWa} Shareshian and Wachs initiated a study of the joint distribution of the Mahonian permutation statistic, major index, and the Eulerian permutation statistic, excedance number.   They showed that the $q$-Eulerian polynomial defined by 
$$\sum_{\sigma \in \S_n} q^{\maj(\sigma)} t^{\exc(\sigma)}$$
can be obtained by taking a certain specialization of the coefficient of $z^n$ in (\ref{hsymeq}) known as the stable principal specialization.
  This enabled them to derive the $q$-analog  of a formula of Euler   that is given in (\ref{expgeneq}).

The work in \cite{ShWa1} led Stanley  to derive a refinement of a formula of Carlitz, Scoville and Vaughan \cite{CaScVa} for the  enumerator of words over the alphabet $\pp$ with no adjacent repeats (see \cite{ShWa}).  These are known as Smirnov words.   In Stanley's refinement the Smirnov words 
are enumerated according to their number of descents and  the enumerator is equal to the symmetric function obtained by applying the standard  involution $\omega$ on the ring of symmetric functions   to   (\ref{hsymeq}).

In another direction, Stanley \cite[Proposition 7.7]{St0} used a formula of Procesi \cite{Pr} to show that (\ref{hsymeq}) is the generating function for the Frobenius characteristic of the representation of the symmetric group $\S_n$ on the cohomology of the toric variety associated with the Coxeter complex of type $A_{n-1}$.

The $q$-Eulerian polynomial, the Smirnov word enumerator, and the   $\S_n$-representation on the cohomology of the toric variety  all have very nice generalizations.  In this paper we discuss these generalizations and explore the relationships among these generalizations and the ramifications of the relationships.   The proofs of the results surveyed here either have appeared in \cite{ShWa} or will appear in \cite{ShWa3}.

The generalization of the $q$-Eulerian numbers,   as defined in \cite{ShWa3}, is obtained 
by considering the joint distribution of  a  Mahonian statistic of Rawlings \cite{Ra} that interpolates between  inversion number and  major index, and a generalized Eulerian statistic of De Mari and Shayman \cite{DeSh}.   We discuss these generalized $q$-Eulerian numbers  and a more general version of them  associated with posets  in Section~\ref{rawsec}.  

The generalization of the Smirnov word enumerator, as defined in \cite{ShWa3}, is a refinement of  Stanley's well known chromatic symmetric function \cite{St3}, which we call a chromatic quasisymmetric function.    In Section~\ref{rawsec} we report on results on
the chromatic quasisymmetric functions, which are proved in \cite{ShWa3}.    In particular we present a refinement of  Gasharov's  \cite{Ga}  expansion of certain chromatic symmetric functions in the Schur basis and a refinement of Chow's \cite{Ch} expansion in the basis of fundamental quasisymmetric functions.   We also present a  conjectured refinement of Stanley's \cite{St3} expansion  in the power-sum basis and a  refinement of the conjecture of Stanley   and Stembridge (see \cite[Conjecture 5.5]{StSte} and \cite[Conjecture 5.1]{St3})  on $e$-positivity.

The generalization of the toric variety associated with the type A Coxeter complex is the regular semisimple Hessenberg variety of type A first studied by De Mari and Shayman \cite{DeSh} and further studied by De Mari, Procesi and Shayman \cite{DePrSh}.     Tymoczko \cite{Ty2} defined a representation of the symmetric group on the cohomology of this Hessenberg variety, which generalizes the representation studied by Procesi and Stanley. This generalization is described in more detail in Section~\ref{hesssec}.

 In Section \ref{cromsec} we present the following extension of the relation between $q$-Eulerian polynomials and Smirnov words.  As noted above, we can associate to a finite poset $P$ a generalized $q$-Eulerian polynomial, which will be denoted by $A_P(q,t)$.  Let $X_G(\x,t)$ denote the chromatic quasisymmetric function of a graph $G$.
 \begin{thm}[Corollary~\ref{maincor}] \label{introth} Let  $G$ be the incomparability graph of a  poset $P$ on $[n]$. Then
$$A_P(q,t)=(q;q)_n \, \Ps( \omega X_G(\x,t)),$$ where
$\Ps$ denotes  stable principal specialization.
\end{thm}
    
     In Section~\ref{hesssec} we discuss the following conjectured generalization of relationship between Smirnov words and the representation of the symmetric group on the cohomology of the toric variety. 
 \begin{con}[Conjecture~\ref{hesschrom}]  \label{introconj} Let $G$ be the incomparability graph of a natural unit interval order $P$.  Then
 $$\omega X_G({\bf x},t)= \sum_{j=0}^{|E(G)|} \ch H^{2j}(\hess(P)) t^j,$$ where $\ch$ denotes the Frobenius characteristic  and $H^{2j}(\hess(P))$ denotes Tymoczko's representation of $\sg_n$ on the degree $2j$ cohomology of the Hessenberg variety $\hess(P)$  associated with $P$.
 \end{con}
This conjecture has the potential of providing solutions to several open problems such as,  
\begin{enumerate}
\item {\it Our  conjecture  that the generalized $q$-Eulerian polynomials are unimodal (Conjecture~\ref{unieucon}).} This would follow from Theorem~\ref{introth} and the hard Lefschetz theorem  applied to Tymoczko's representation on the cohomology of the Hessenberg variety. 
\item {\it Tymoczko's   problem  of finding  a decomposition of her representation into irreducibles \cite{Ty2}.}  Such a decomposition would be provided by the expansion of  the chromatic quasisymmetric functions in the Schur function basis given in Theorem~\ref{schurcon}.  
\item {\it Stanley and Stembridge's well-known conjecture  that the chromatic symmetric functions (for  unit interval orders) are $e$-positive \cite{StSte,St3}.}  This would be equivalent to Conjecture~\ref{tympermcon} that  Tymoczko's representation is a permutation representation in which each point stabilizer is a Young subgroup.   
\end{enumerate}

\subsection{Preliminaries} \label{prelimsec} We assume that the reader is familiar with some basic notions from combinatorics and symmetric function theory.  All terms used but not defined in this paper are defined in \cite{St1} or \cite{St2}.   

Let $\sigma \in \sg_n$, where  $\S_n$ denotes the symmetric group on $[n]:=\{1,2,\dots,n\}$. The 
excedance number of $\sigma$ is given by
$$\exc(\sigma) :=|\{i \in [n-1] : \sigma(i) > i\}|.$$   The descent set of $\sigma$  is
given by $$\Des(\sigma) :=  \{i \in [n-1] : \sigma(i) > \sigma(i+1)\}$$ and the descent number and 
major index are $$\des(\sigma) := |\Des(\sigma) |\,\, \mbox{  and }\,\,
 \maj(\sigma):= \sum_{i\in \Des(\sigma)} i .$$  It is well-known that the permutation statistics $
 \exc$ and $\des$ are   equidistributed on $\S_n$.  The common generating functions for these statistics are called Eulerian polynomials.  That is, the Eulerian polynomials are defined as $$A_{n}(t) :=\sum_{\sigma \in \S_n}t^{\des(\sigma)} = \sum_{\sigma \in \S_n}t^{\exc(\s)} .$$
 Any permutation statistic with generating function $A_n(t)$ is called an Eulerian statistic. The Eulerian numbers are the coefficients of the Eulerian polynomials; they are defined as
 $$a_{n,j} := |\{\sigma \in \S_n : \des(\sigma) = j\}| =| \{\sigma \in \S_n : \exc(\sigma) = j\}|,$$
 for $0 \le j\le n-1$.

It is also well-known that the  major index is equidistributed with the inversion index defined as
$$\inv(\sigma) = |\{(i,j) \in [n] \times [n] : i < j \mbox{ and } \sigma(i) > \sigma(j)\}|$$ and that
$$\sum_{\sigma \in \S_n}q^{\maj(\sigma)} = \sum_{\sigma \in \S_n}q^{\inv(\s)} = [n]_q!,$$
where 
$$[n]_q := 1+q+ \dots + q^{n-1} \mbox{ and } [n]_q! := [n]_q [n-1]_q \dots [1]_q. $$
Any permutation statistic equidistributed with $\maj $ and $\inv$ is called a Mahonian permutation statistic.  

 Let  $R$ be a    ring  with a partial order relation, e.g., $\Q$ with its usual total order.  Given  a sequence $(a_0,a_1,\dots,a_n)$ of elements of $R$ we say that the sequence  is  {\em palindromic} with center of symmetry $\frac n 2$ if
$a_j = a_{n-j}$  for $0 \le j \le n$.   The sequence is  said to be {\em  unimodal}  if 
$$a_0 \le a_1\le \cdots \le a_c \ge a_{c+1} \ge a_{c+2} \ge \cdots \ge a_n$$
for some $c$.  The sequence is said to be {\em positive} if $a_i \ge 0$ for each $i$.   We say that the polynomial
$A(t):=a_0 + a_1 t +\cdots  + a_{n} t^{n}\in R[t]-\{0\}$ is positive, palindromic and unimodal with center of symmetry $\frac {n}2$ if   $(a_0,a_1,\dots,a_{n})$ has these properties.

Let $\Lambda_\Q$ be the ring of symmetric functions in variables $x_1,x_2,\dots$ with coefficients in $\Q$ and  let $b$ be any basis for  $\Lambda_\Q$.  If $f$ is a symmetric function whose expansion in $b$ has nonnegative coefficients, we say that $f$ is $b$-positive. This induces a partial order relation on the ring $\Lambda_\Q$  given by  $f \le_b g $ if $g-f$ is $b$-positive.   We  use this partial order  to define  $b$-unimodality of sequences of symmetric functions in $\Lambda_\Q$ and polynomials in $\Lambda_\Q{[t]}$.   Similarly, when we say that a sequence of polynomials in $\Q[q_1,\dots,q_m]$ or a polynomial in $\Q[q_1,\dots, q_m][t]$ is unimodal, we are using the partial order relation on $\Q[q_1,\dots,q_m]$ defined by 
 $f({\bf q}) \le g({ \bf q})$   if $g({\bf q}) -f({\bf q}) $ has nonnegative coefficients.  
 
 The bases for $\Lambda_\Q$ that are relevant here are the 
 Schur basis $(s_\lambda)_{\lambda \in \Par}$, the elementary symmetric function basis
$(e_\lambda)_{\lambda \in \Par}$, the complete homogeneous symmetric function basis $(h_\lambda)_{\lambda \in \Par}$, and the power-sum symmetric function basis $(z_\lambda^{-1}p_\lambda)_{\lambda \in \Par}$,
 where $\Par$ is the set of all integer partitions, $$z_\lambda:=\prod_i m_i(\lambda)! \, i^{m_i(\lambda)}$$ and $m_i(\lambda)$ is the number of parts of $\lambda$ equal to $i$.    These are all integral bases in that every symmetric function in the ring $\Lambda_\Z$ of symmetric functions over $\Z$ is an integral linear combination of the basis elements.

The next two propositions give basic tools for establishing unimodality.
\begin{prop}[see {\cite[Proposition 1]{St0}}] \label{tooluni} Let $A(t)$ and $B(t)$ be positive, palindromic and unimodal polynomials in $\Q[t]$ with respective centers of symmetry $c_A$ and $c_B$.  Then
\begin{enumerate}
\item $A(t) B(t) $ is positive, palindromic and unimodal with center of symmetry $c_A +c_B$.
\item If $c_A=c_B$ then $A(t) +B(t)$ is positive, palindromic and unimodal with center of symmetry $c_A$.
\end{enumerate}
\end{prop}

\begin{prop} \label{tooluni2}
Let  $b:= (b_\lambda)_{\lambda}$ be a basis for $\Lambda_\Q
$.  The polynomial $f(t) \in \Lambda_{\Q}[t]$ is $b$-positive, palindromic and $b$-unimodal with center of symmetry $c$ if and 
only if each coefficient $a_\lambda(t)$ in the expansion $f(t) = \sum_\lambda a_\lambda(t) b_\lambda$ is 
positive, palindromic and unimodal with center of symmetry $c$.
\end{prop}

For example, from Proposition~\ref{tooluni} (1) we see that $[n]_t!$ is positive, palindromic and unimodal with center of symmetry $\frac {n(n+1)}{4}$.  From Propositions~\ref{tooluni} and~\ref{tooluni2} we have that $[5]_t [2]_t h_{(5,2)} + [4]_t [3]_t h_{(4,3)}$ is $h$-positive, palindromic and $h$-unimodal with center of symmetry $5/ 2$.

 Recall that  the Frobenius characteristic  $\ch$  is a ring isomorphism from the ring of virtual representations  of the symmetric groups to $\Lambda_\Z$.  It takes the irreducible Specht module $S^{\lambda}$ to the Schur function $s_\lambda$.  
 As is customary, we use  $\omega$ to denote the involution on $\Lambda_\Q$ that takes $h_i$ to $e_i$.  Recall $\omega(s_\lambda) = s_{\lambda^\prime}$ for all $\lambda \in \Par$,  where  $\lambda^\prime$ denotes the conjugate shape. 
  
For any ring $R$, let $\cq_R$ be the ring of quasisymmetric functions in variables $x_1,x_2,\ldots$ with coefficients in  $R$.  
For $n \in \pp$ and $S \subseteq [n-1]$, let $D(S)$ be the set of all functions $f:[n] \rightarrow \pp$ such that
\begin{itemize}
\item $f(i) \geq f(i+1)$ for all $i \in [n-1]$, and \item $f(i)>f(i+1)$ for all $i \in S$.
\end{itemize}
Recall Gessel's  {\em fundamental quasisymmetric function}
\[
F_{n,S}:=\sum_{f \in D(S)}\x_f,
\]
where $\x_f := x_{f(1)}x_{f(2)} \cdots x_{f(n)}$.
It is straightforward to confirm that $F_{n,S} \in \cq_{\Z}$.  In fact (see \cite[Proposition 7.19.1]{St2}),  $\{F_{n,S}:S \subseteq [n-1]\}$ is a basis for the free $\Z$-module of homogeneous degree $n$ quasisymmetric functions with coefficients in $\Z$.

 See \cite{St2} for further information on symmetric functions, quasisymmetric functions and representations.

\section{$q$-Eulerian numbers and toric varieties} \label{eultorsec}
It is well known that the Eulerian numbers $(a_{n,0},a_{n,1},\dots,a_{n,n-1})$  form the $h$-vector  of the Coxeter complex $\Delta_n$ of the symmetric group $\sym_n$.     Danilov and Jurciewicz (see \cite[eq.~(26)]{St0}) showed that the $j$th entry of the h-vector of any convex simplicial polytope $P$  is equal to  $\dim H^{2j}(\v(P))$ where $\v(P)$ is the toric variety associated with $P$ and $H^{i}(\v(P))$ is the degree $i$ singular cohomology of $\v(P)$ over $\C$.  (In odd degrees cohomology  vanishes.)   It follows that  
\begin{equation} \label {eultoreq} a_{n,j} = \dim H^{2j}(\v(\Delta_n)) ,\end{equation}
 for all $j = 0,\dots,n-1$.
 
 It is also well known that the Eulerian numbers  are palindromic and unimodal.
 Although there are elementary ways to prove this result, palindromicity and unimodality  can be viewed as   consequences of the hard Lefschetz theorem applied to the toric variety 
 $\v(\Delta_n)$.   
 
  The classical hard Lefschetz theorem states that if $\v$ is a projective  smooth variety  of complex dimension $n$ then there is a linear operator
$\phi$ on $H^*(\v)$ sending each $H^{i}(\v)$ to $H^{i+2}(\v)$ such that  for $0\le i \le 2n-2$, the map $\phi^{n-i}: H^{i}(\v) \to H^{2n-i}(\v)$ is a bijection.  It follows from this that the sequence of even degree Betti numbers and the sequence of odd degree Betti numbers are palindromic and unimodal.  Since $\v(\Delta_n)$ is a smooth variety, it follows that the Eulerian numbers form a palindromic and unimodal sequence.  See \cite{St-1,St0} for further information on the use of the hard Lefschetz theorem in combinatorics.

\subsection{Action of the symmetric group} \label{torsec} The symmetric group $\sym_n$ acts naturally on the Coxeter complex $\Delta_n$ and this induces a representation of $\sym_n$ on cohomology of the toric variety $\v(\Delta_n)$.  
Procesi \cite{Pr} derived a recursive formula for this representation, which Stanley \cite[see Proposition 12]{St0} used to obtain the following formula for the Frobenius characteristic of the representation, 
\begin{equation} \label{symeuler}1+ \sum_{n\ge 1} \sum_{j= 0} ^{n-1}\ch H^{2j}(\v(\Delta_n)) t^j z^n = {(1-t) H(z) \over H(tz) - t H(z)} .\end{equation}  
(See   \cite{Ste1,Ste2,DoLu, Le} for related work on this representation.)     This is a symmetric function analog of the classical formula of Euler,
\begin{equation} \label{expgeneqn}
1+ \sum_{n\ge 1} \sum_{j= 0} ^{n-1}a_{n,j}t^j\frac{z^n}{n!}=\frac{(1-t)e^z}{e^{tz}-te^z}.
\end{equation}

An immediate consequence of (\ref{symeuler})  is 
that $\ch H^{2j}(\v(\Delta_n))$ is $h$-positive, which is equivalent to saying  that the linear representation of $\S_n$ on $H^{2j}(\v(\Delta_n))$ is obtained from an action on some set in which each point stabilizer is a Young subgroup. 
Formula (\ref{symeuler}) has a number of additional interesting, but not so obvious, consequences. We discuss some of them below.  

Stembridge \cite{Ste1}  uses (\ref{symeuler})  to characterize the multiplicity of each irreducible $\sym_n$-module in $H^{2j}(\v(\Delta_n))$ in terms of marked  tableaux.   A {\em marked tableau} of shape $\lambda$ is a semistandard tableau of shape $\lambda$ with entries in $\{0,1,\dots,k\}$ for some $k \in \N$, such that each $a \in [k]$ occurs at least twice in the tableau, together with a mark on one occurrence of each $a \in [k]$; the marked occurrence cannot be the leftmost occurrence.  The index of a marked $a$ is the number of occurrences of $a$ to the left of the marked $a$.  The index of a marked tableau is the sum of the indices of the marked entries of the tableau.  

\begin{thm}[Stembridge \cite{Ste1}] \label{stemdecomp} Let $\lambda \vdash n$.  The multiplicity of the irreducible Specht module $S^\lambda$ in $H^{2j}(\v(\Delta_n)) $  is the number of marked tableaux of shape $\lambda$ and index $j$.
\end{thm}

Recently,  Shareshian and Wachs \cite{ShWa3} obtained  a  different formula for the multiplicity of $S^\lambda$  in terms of a different type of tableau (see Corollary~\ref{swdecomp}).

The character of $H^{2j}(\v(\Delta_n))$ can be obtained from the following expansion in the $p$-basis.

\begin{thm}[Stembridge \cite{Ste1}, Dolgachev and Lunts \cite{DoLu}] \label{pdecomp} For all $n\ge 1$, $$\sum_{j=0}^{n-1} \ch H^{2j}(\v(\Delta_n)) t^j = \sum_{\lambda \vdash n} A_{l(\lambda)}(t) \prod_i [\lambda_i]_t \,\, z_\lambda^{-1} p_\lambda\,.$$
\end{thm}

\subsection{Expansion in the basis of fundamental quasisymmetric functions}

 In their study of $q$-Eulerian numbers \cite{ShWa1,ShWa}, Shareshian and Wachs    obtain a formula for the expansion of  the right hand side of (\ref{symeuler}) in the basis of fundamental quasisymmetric functions. In order to describe this 
decomposition  the notion of $\Dex$ set of a permutation is needed.  For $n \ge 1$, we set \[
[\ov{n}]:=\{\ov{1},\ldots,\ov{n}\}
\]
and totally order the alphabet $[n] \cup [\ov{n}]$ by
\begin{equation} \label{order1}
\ov{1}<\ldots<\ov{n}<1<\ldots<n.
\end{equation}
For a permutation $\sigma=\sigma_1\ldots\sigma_n \in \S_n$, we
define $\ov{\sigma}$ to be the word over  alphabet $[n] \cup
[\ov{n}]$ obtained from $\sigma$ by replacing $\sigma_i$ with
$\ov{\sigma_i}$ whenever $i \in \Exc(\sigma)$.  For example, if
$\sigma={\rm {531462}}$ then $\ov{\sigma}={\rm
{\ov{5}\ov{3}14\ov{6}2}}$.  We define a descent in a word
$w=w_1\ldots w_n$ over any totally ordered alphabet to be any $i
\in [n-1]$ such that $w_i>w_{i+1}$ and let $\Des(w)$ be the set of
all descents of $w$.  Now, for $\sigma \in \S_n$, we define
\[
\Dex(\sigma):=\Des(\ov{\sigma}).
\]
For example, $\Dex({\rm {531462}})=\Des({\rm {\ov{5}\ov{3}14\ov{6}2}})=\{1,4\}$.

\begin{thm}[Shareshian and Wachs {\cite[Theorem 1.2]{ShWa}}] \begin{equation} \label{eulfund}  1+ \sum_{n\ge 1} \sum_{\sigma \in \S_n} F_{n,\Dex(\sigma)} t^{\exc(\sigma)} z^n = {(1-t) H(z) \over H(tz) - t H(z)}  .\end{equation}  
\end{thm}

\begin{cor}[{\cite[Theorem 7.4]{ShWa}}] \label{funddecompth} For  $0\le j\le n-1$,
$$\ch H^{2j}(\v(\Delta_n)) = \sum_{\scriptsize\begin{array}{c} \sigma \in \sym_n \\ \exc(\sigma) = j\end{array}} F_{n,\Dex(\sigma)} \,.$$
 \end{cor}
 
Expansion in the basis of fundamental quasisymmetric functions is useful because it can yield interesting results about permutation statistics via specialization. 
For a quasisymmetric function $Q(x_1,x_2,\ldots) \in \cq_R$, the stable principal specialization $\Ps(Q) \in R[[q]]$ is, by definition, obtained from $Q$ by substituting $q^{i-1}$ for $x_i$ for each $i \in \pp$. Now (see for example \cite[Lemma 7.19.10]{St2}), for any $n \in \pp$ and $S \subseteq [n-1]$, we have
\begin{equation} \label{sps}
\Ps(F_{n,S})=(q;q)_n^{-1} \,\,q^{\sum_{i \in S}i},
\end{equation}
where $$(p;q)_n = \prod_{j=1}^n (1-pq^{j-1}).$$
From Corollary~\ref{funddecompth}  and the fact that that $\sum_{i \in \Dex(\sigma) } i = \maj(\sigma) - \exc(\sigma)$ (see \cite[Lemma 2.2]{ShWa}) we obtain the following q-analog of (\ref{eultoreq}).
\begin{cor} \label{pscor} For $0 \le j \le n-1$, $$ (q;q)_{n} \,\,\Ps( \ch H^{2j}(\v(\Delta_n)) )=a_{n,j}(q) ,$$ where  \begin{equation} \label{qeulerdef1} a_{n,j}(q):= \sum_{\scriptsize\begin{array}{c} \sigma \in \sym_n \\ \exc(\sigma) = j\end{array}} q^{\maj(\sigma)-\exc(\sigma)}.\end{equation}\end{cor}
The $q$-Eulerian numbers $a_{n,j}(q)$ defined in Corollary~\ref{pscor}  were initially studied in \cite{ShWa1,ShWa} and have been further studied in \cite{ShWa2.5, FoHa1, FoHa2, SaShWa, HeWa}.  
Define the $q$-Eulerian polynomials,
 \begin{equation} \label{qeulerdef2}A_n(q,t) := \sum_{j=1}^{n-1} a_{n,j}(q) t^j = \sum_{\sigma \in \sym_n} q^{\maj(\sigma) -\exc(\sigma)} t^{\exc(\sigma)}.\end{equation}
 In \cite{ShWa1,ShWa}, the following $q$-analog of Euler's formula (\ref{expgeneqn}) was obtained via specialization of (\ref{eulfund}),
 \begin{equation} \label{expgeneq}
1+\sum_{n \geq 1}A_n(q,t)\frac{z^n}{[n]_q!}=\frac{(1-t)\exp_q(z)}{\exp_q(tz)-t\exp_q(z)},
\end{equation}
where $$\exp_q(z) := \sum_{n \ge 0} \frac {z^n}{[n]_q!}.$$

It is surprising that when one evaluates $a_{n,j}(q)$ at  any $n$th root of unity, one always gets a positive integer.

\begin{thm}[Sagan, Shareshian, and Wachs {\cite[Corollary 6.2]{SaShWa}}]  \label{unity} Let $dm=n$ and let $\xi_d$ be any primitive $d$th root of unity.   Then
$$A_n(\xi_d,t) = A_m(t) [d]^m_t.$$  Consequently $A_n(\xi_d,t)$ is a positive, palindromic, unimodal polynomial in $\Z[t]$.
\end{thm}

Theorem~\ref{pdecomp}, Corollary~\ref{pscor} and the following result, which is  implicit in \cite{De} and  stated explicitly in \cite{SaShWa}, can be used to give a proof of Theorem~\ref{unity}.

\begin{lemma} \label{thdes} Suppose $u(q) \in \Z[q]$ and there exists  a homogeneous symmetric function $U$ of degree $n$ with coefficients in $\Z$ such that   
 $$u(q) =(q;q)_n\,\, \Ps (U).$$ If $dm =n$  then $u(\xi_d)$ is the coefficient of $z_{d^m}^{-1} p_{d^m}$ in the expansion of $U$ in the power-sum basis.  
\end{lemma}

\subsection{Unimodality} 
Let us now consider what the hard Lefschetz theorem tells us about
the sequence $( \ch H^{2j}(\v(\Delta_n)))_{0\le j\le n-1}$.  Since the action of the symmetric group $\sym_n$ commutes with the hard Lefschetz map $\phi$ (see \cite[p. 528]{St0}), we can conclude 
that for $0 \le i \le n$ the map $\phi:H^{i}(\v(\Delta_n)) \to H^{i+2}(\v(\Delta_n))$ is an $\sym_n$-module monomorphism.  Hence, by Schur's lemma, for each $\lambda \vdash n$, the multiplicity of the Specht module  $S^\lambda$ in $H^{i}(\v(\Delta_n))$ is less than or equal to the multiplicity in $H^{i+2}(\v(\Delta_n))$.  Equivalently, the coefficient of the Schur function $s_\lambda$ in the expansion of $\ch H^{i}(\v(\Delta_n))$ in the Schur basis is less than or equal to the coefficient in $\ch H^{i+2}(\v(\Delta_n))$.   
Hence, it follows from the hard Lefschetz theorem that the sequence $( \ch H^{2j}(\v(\Delta_n)))_{0\le j\le n-1}$ is palindromic and Schur-unimodal.

The following lemma is useful for establishing unimodality of polynomials in $\Q[q][t]$.
\begin{lemma}[see for example {\cite[Lemma 5.2]{ShWa}}]  \label{hposqpos}  If $U$ is a Schur-positive homogeneous symmetric function of degree $n$  then $$(q;q)_{n} \Ps (U)$$ is a polynomial in $q$ with nonnegative coefficients.
 \end{lemma}
 
 It follows from Lemma~\ref{hposqpos} and Corollary~\ref{pscor} that the palindromicity and Schur-unimodality of $( \ch H^{2j}(\v(\Delta_n)))_{0\le j\le n-1}$ can be specialized to yield the following $q$-analog of the palindromicity and unimodality of the Eulerian numbers.  

\begin{thm}[Shareshian and Wachs {\cite[Theorem 5.3]{ShWa}}] \label{unieuth} For any $n \in \pp$, the sequence $(a_{n,j}(q))_{0\le j\le n-1}$ is palindromic and unimodal.  
\end{thm}

We remark that it is not really necessary to use the hard Lefschetz theorem to establish the above mentioned unimodality results. 
In fact, by manipulating the right hand side of (\ref{symeuler}) we obtain (\ref{formQ}) 
below\footnote{This formula is  different from a similar looking formula given in \cite[Corollary 4.2] {ShWa}.}. The  palindromicity and h-unimodality consequences follow by Propositions~\ref{tooluni} and \ref{tooluni2}.  

\begin{thm}[Shareshian and Wachs \cite{ShWa3}]  \label{hunimodal}  For all $n \ge 1$, 
\begin{equation} \label{formQ}  \sum_{j=0}^{n-1} \ch H^{2j}(\v(\Delta_n)) t^j = \sum_{m=1}^{\lfloor {n+1 \over 2} \rfloor} t^{m-1}  \sum_{\substack{
k_1,\dots, k_m \ge 2 \\ \sum k_i = n+1}}\,\,\, \prod_{i=1}^m  [k_i-1]_{t} \,\,h_{k_1-1}h_{k_2} \cdots h_{k_m}\end{equation}
Consequently, the sequence $( \ch H^{2j}(\v(\Delta_n)))_{0\le j\le n-1}$ is $h$-positive, palindromic and $h$-unimodal.
\end{thm}

By specializing (\ref{formQ}) we obtain the following result.  Recall that the $q$-analog of the multinomial coefficients is defined by $$ \left[\begin{array}{c} n \\k_1, \dots, k_m\end{array}\right]_q\, := \frac {[n]_q!}{ [k_1]_q! \cdots [k_m]_q!}$$ for all $k_1,\dots,k_m \in \N$ such that $\sum_{i=1}^n k_i = n$.
\begin{cor}For all $n \ge 1$, $$
 A_n(q,t) =  \sum_{m=1}^{\lfloor {n+1 \over 2} \rfloor} t^{m-1}  \sum_{\substack{
k_1,\dots, k_m \ge 2 \\ \sum k_i = n+1}}\,\,\, \prod_{i=1}^m  [k_i-1]_{t} \,\, \left[\begin{array}{c} n \\ k_1-1,k_2, \dots, k_m\end{array}\right]_q\, .
$$
Consequently, $A_n(q,t)$ is a palindromic and unimodal polynomial in $t$.
\end{cor}

\subsection{Smirnov words} \label{smir} The expression on the right hand side of (\ref{symeuler}) has appeared in several other contexts (see \cite[Section 7] {ShWa}).  We mention just  one of these contexts here.

A Smirnov word is a word  with no equal adjacent letters.   Let $W_n$ be the set of all Smirnov words of length $n$ over alphabet $\PP$. Define the enumerator
$$Y_{n,j}(x_1,x_2,\dots):= \sum_{\scriptsize \begin{array}{c} w\in  W_n \\ \des(w) = j\end{array}} {\bf x}_w,$$
 where
${\bf x}_w:=x_{w(1)} \cdots x_{w(n)}$. (Here we are calculating $\des(w)$ using the standard total order on $\pp$.)
A formula for the generating function of $Y_n:=\sum_{j=0}^{n-1} Y_{n,j}$ was initially obtained by  Carlitz, Scoville and Vaughan  \cite{CaScVa} and further studied by Dollhopf, Goulden and Greene  \cite{DoGoGr} and Stanley \cite{St3}.  Stanley pointed out to us  the following refinement of this formula.   This refinement follows from results in \cite[Section 3.3]{ShWa} by P-partition reciprocity.

\begin{thm}[Stanley (see {\cite[Theorem  7.2]{ShWa}})] \label{stanth} 

\begin{equation}\label{carlg} \sum_{n,j \ge 0} Y_{n,j}({\bf x})t^j z^n = {(1-t) E(z) \over E(zt) - t E(z)} ,\end{equation} where  $E(z) := \sum_{n \ge 0} e_n z^n$ and $e_n$ is the elementary symmetric function of degree $n$.
\end{thm}

\begin{cor}\label{stancor} For $0 \le j \le n-1$,
$$ \ch H^{2j}(\v(\Delta_n)) = \omega Y_{n,j}.$$
\end{cor}

\section{Rawlings major index and  generalized Eulerian numbers} \label{rawsec}

In \cite{Ra}, Rawlings studies Mahonian permutation statistics that interpolate between the major index and the inversion index.   
  Fix $n \in \pp$ and $k \in [n]$.  For $\sigma \in \S_n$, set
\begin{eqnarray*}
\Des_{\geq k}(\sigma)&:=&\{i \in [n-1]:\sigma(i)-\sigma(i+1) \geq k\}, \\
\maj_{\ge k} (\sigma)& := & \sum_{i \in \Des_{\ge k}} i ,
\end{eqnarray*}

\begin{eqnarray*}
\Inv_{<k}(\sigma)&:=&\{(i,j) \in [n] \times [n]:i<j \mbox{ and } 0<\sigma(i)-\sigma(j)<k\} ,\\
\inv_{<k}(\sigma)&:=& |\Inv_{<k}(\sigma)|.
\end{eqnarray*}

Now define the {\it Rawlings major index} to be
\[
\rmaj_k(\sigma):=\inv_{<k}(\sigma)+\maj_{\ge k}(\sigma).
\]
Note that $\rmaj_1$ and $\rmaj_n$ are, respectively, the well studied major index $\maj$ and inversion number $\inv$.  Rawlings shows in \cite{Ra} that each $\rmaj_k$ is a Mahonian statistic, that is,
\bq \label{raweq}
\sum_{\sigma \in \S_n}q^{\rmaj_k(\sigma)}=[n]_q!\,.
\eq
For a proof of (\ref{raweq}) that is different from that of Rawlings and a generalization of (\ref{raweq}) from permutations  to labeled trees, see \cite{LiWa}.

In  \cite{DeSh}  De Mari and Shayman introduce a class of numbers, closely related to the Rawlings major index, which they call generalized Eulerian numbers.    For $k \in [n]$  and $0\le j $ define the De Mari-Shayman generalized Eulerian numbers to be $$a_{n,j}^{(k)} := |\{\sigma \in \S_n : \inv_{< k} (\sigma)  = j\}|.$$
Note that  $\inv_{<2}(\sigma) = \des(\sigma^{-1})$, which implies that $a_{n,j}^{(2)} = a_{n,j}$,  justifying the name ``generalized Eulerian numbers".    De Mari and Shayman introduce the generalized Eulerian numbers in connection with their study of Hessenberg varieties.  It follows from their work and the hard Lefschetz theorem that for each fixed $k\in [n]$, the generalized Eulerian numbers $(a^{(k)}_{n,j})_{j \ge 0}$ form a palindromic unimodal sequence of numbers.

In \cite{ShWa3} Shareshian and Wachs consider  a $q$-analog of the De Mari-Shayman generalized Eulerian numbers defined  for $k \in [n]$  and $0\le j $ by
$$a_{n,j}^{(k)}(q) := \sum_{\scriptsize{\begin{array}{c} \sigma \in S_n \\ \inv_{< k} (\sigma) = j \end{array}} } q^{\maj_{\ge k} (\sigma)}.$$
Similarly, the generalized $q$-Eulerian polynomials are defined by
$$A^{(k)}_n(q,t) = \sum_{\sigma \in \S_n} q^{\maj_{\ge k} (\sigma)} t^{ \inv_{< k} (\sigma) }.$$
Now (\ref{raweq}) is equivalent to,
$$A_n^{(k)}(q,q) = [n]_q!.$$

We now consider the question of whether unimodality and other known properties of the generalized Eulerian numbers $a_{n,j}^{(k)}$ extend to the generalized $q$-Eulerian numbers $a^{(k)}_{n,j}(q)$.   Although it is not at all obvious, it turns out that  in the case $k=2$,  the generalized $q$-Eulerian numbers   are equal to the $q$-Eulerian numbers  defined in  (\ref{qeulerdef1}),
$$a^{(2)}_{n,j}(q) = a_{n,j}(q)$$
(see Theorem~\ref{eqEu}). Hence by Theorem~\ref{unieuth}, unimodality holds when $k=2$.    We  consider this question in a more general setting, which we now describe.

Let $E(G)$ denote the edge set of a graph $G$.   For $\sigma \in \S_n$ and $G$ a graph with vertex set $[n]$, the {\em $G$-inversion set} of $\sigma$ is $$ \Inv_G(\sigma):=\{ (i,j) : i<j, \,\,\sigma(i) > \sigma(j) \mbox{ and } \{\sigma(i),\sigma(j) \} \in E(G) \} $$ 
and the {\em $G$-inversion number} is $$\inv_G(\sigma) := | \Inv_G(\sigma)|.$$
For $\sigma \in \S_n$ and $P$ a poset on $[n]$, the {\em $P$-descent set} of $\sigma$ is 
$$ \Des_P(\sigma ):= \{ i \in [n-1] : \sigma(i) >_P \sigma(i+1) \},$$ and the {\em $P$-major index } is $$\maj_P(\sigma)  := \sum_{i \in \Des_P(\sigma)} i.$$

 Define the {\em incomparability graph} $\inc(P)$  of a poset $P$ on $[n]$ to be  the graph with vertex set $[n]$ and edge set  $\{\{a,b\} : a \not\le_P b \mbox{ and } b \not\le_P a \}$.
 For  $0 \le j \le |E(\inc(P))|$ define the $( q,P)$-Eulerian numbers,
$$a_{P,j}(q) := \sum_{\scriptsize\begin{array}{c} \sigma \in \S_n \\ \inv_{\inc(P)}(\sigma) = j\end{array}} q^{\maj_P(\sigma)},$$
and  the  $(q,P)$-Eulerian polynomials
$$A_P(q,t) := \sum_{j=0}^{|E(\inc(P))|} a_{P,j}(q) t^j = \sum_{\sigma \in \S_n} q^{\maj_P(\sigma)} t^{\inv_{\inc(P)}(\sigma)}.$$

Define $P_{n,k}$ to be the poset on vertex set $[n]$ such that $i<j$ in $P_{n,k}$ if and only if $j-i \geq k$ and let $G_{n,k}$ be the incomparability graph of $P_{n,k}$.  Then 
\[
\Des_{\geq k}(\sigma)=\Des_{P_{n,k}}(\sigma)
\]
and
\[
\Inv_{<k}(\sigma)=\Inv_{G_{n,k}}(\sigma).
\]
Hence if $P=P_{n,k} $, we have $a_{P,j}(q) = a_{n,j}^{(k)}(q)$ and $A_P(q,t) = A^{(k)}_{n}(q.t)$.

The unimodality property that holds for the $( q,P)$-Eulerian numbers  in case $P =P_{n,2}$ does not hold for general $P$.  Indeed, consider the poset $P$ on $[3]$ whose only relation is $1<2$.   We have
\begin{equation} \label{nonnatint} A_P(q,t) := (1+q) + 2t +(1+q^2) t^2,\end{equation}
which is neither palindromic as a polynomial in $t$, nor unimodal.  The property $A_P(q,q) = [n]_q! $ fails as well.
However there is a very nice class of posets for which unimodality seems to hold.  

A {\em  unit interval order} is a  poset  that is isomorphic to a finite collection $\ci$ of  intervals $[a,a+1]$ on the real line, partially ordered by the relation $[a,a+1]<_\ci[b,b+1]$ if $a+1<b$.  Define a {\em natural unit interval order} to be a poset $P$ on   $[n]$ that satisfies the following conditions
\begin{enumerate} 
\item $x<_P y$ implies $x<y$ in the natural order on $[n]$
\item if the direct sum $\{x<_P z\} + \{y\}$ is an induced subposet of $P$ then  $x<y<z$ in the natural order on $[n]$.
\end{enumerate}
 It is not difficult to see that  every natural unit interval order is a unit interval order and that every unit interval order is isomorphic to a unique natural unit interval order. 
 The poset $P_{n,k}$ is an example of a natural unit interval order.

It is a consequence of  a result of Kasraoui  \cite[Theorem 1.8]{Ka} that if $P$ is a natural unit interval order  then $\inv_{\inc(P)}+\maj_{P}$ is Mahonian.   In other words, if $P$ is a natural unit interval order  on $[n]$ then
\[
A_P(q,q) = [n]_q! \,\,.
\]

From the work of De Mari, Procesi and Shayman \cite{DePrSh} we have the following result.  The proof  is discussed in Section~\ref{hesssec}.
\begin{thm}[see {\cite[Exercise 1.50 (f)]{St1}}] \label{hesseuth} Let $P$ be a natural unit interval order.  Then the $P$-Eulerian polynomial  $A_P(1,t)$  is  palindromic and unimodal.
\end{thm}

Palindromicity of the $q$-analog $A_P(q,t)$ follows from results discussed in Section~\ref{cromsec} and unimodality  is implied by a conjecture discussed in Section~\ref{cromsec}.

\begin{thm}[Shareshian and Wachs \cite{ShWa3}] \label{palinAth}
Let $P$ be a natural unit interval order on $[n]$.  Then the sequence $(a_{P,j}(q))_{0\le j \le |E(\inc(P))|}$ is  palindromic.
\end{thm}

\begin{con} \label{unieucon} Let $P$ be a natural unit interval order on $[n]$.  Then the palindromic sequence $(a_{P,j}(q))_{0\le j \le |E(\inc(P))|}$ is   unimodal.
\end{con}

The following conjecture  generalizes Theorem~\ref{unity}.

\begin{con} \label{genunity}Let $P$ be a natural unit interval order on $ [n]$.  If $dm = n$ then there is a polynomial $B_{P,d}(t) \in \N[t]$ such that $$A_P(\xi_d,t) = B_{P,d}(t) [d]_t^m.$$ Moreover,  $(a_{P,j}(\xi_d))_{0\le j \le |E(\inc(P))|}$ is a palindromic unimodal sequence of positive integers.
\end{con}

  In addition to the case $P = P_{n,2}$ (Theorems~\ref{unieuth} and \ref{unity}), these conjectures have been verified for $P_{n,k}$ when $k=1,n-2,n-1,n$, and by computer for all $k$ when $n \le 8$.  An approach to proving them in general will be presented in Sections~\ref{cromsec} and \ref{hesssec}.

\section{Chromatic quasisymmetric functions} \label{cromsec}
\subsection{Stanley's chromatic symmetric functions} Let $G$ be a graph with vertex set $[n]$ and edge set $E=E(G) \subseteq {{[n]} \choose {2}}$.  A {\em proper $\pp$-coloring} of $G$ is a function $c$ from $[n]$ to the set $\pp$ of positive integers such that whenever $\{i,j\} \in E$ we have $c(i) \neq c(j)$.  
Given any function $c:[n] \rightarrow \pp$, set
\[
\x_c:=\prod_{i=1}^{n}x_{c(i)}.
\]
Let $C(G)$ be the set of proper $\pp$-colorings of $G$.  In \cite{St3}, Stanley defined the {\em chromatic symmetric function} of $G$,
\[
X_G(\x):=\sum_{c \in C(G)}\x_c.
\]
It is straightforward to confirm that $X_G \in \Lambda_{\Z}$.  The chromatic symmetric function is a generalization of the chromatic polynomial $\chi_G:\pp \to \pp$, where $\chi_G(n)$ is the number of proper colorings of $G$ with $n$ colors.   Indeed,  $X_G(1^n) = \chi_G(n)$, where $X_G(1^n)$ is the specialization of $X_G(\x)$   obtained by setting $x_i=1$ for $1\le i \le n$ and $x_i=0$ for $i > n$.   Chromatic symmetric functions are studied in various papers, including \cite{St3,St4,Ga,Ch,MaMoWa}.

We recall Stanley's description of the power sum decomposition of $X_G(\x)$, for arbitrary $G$ with vertex set $[n]$.  We call a partition $\pi=\pi_1|\ldots|\pi_l$ of $[n]$ into nonempty subsets {\em $G$-connected} if the subgraph of $G$ induced on each block $\pi_i$ of $\pi$ is connected.  The set $\Pi_{G,n}$ of all $G$-connected partitions of $[n]$ is partially ordered by refinement (that is, $\pi \leq \theta$ if each block of $\pi$ is contained in some block of $\theta$).  We write $\mu_G$ and $\hat{0}$, respectively, for the M\"obius function on $\Pi_{G,n}$, and the minimum element $1|\ldots|n$ of $\Pi_{G,n}$.  For $\pi \in \Pi_{G,n}$, we write $\pa(\pi)$ for the partition of $n$ whose parts are the sizes of the blocks of $\pi$. 

\begin{thm}[Stanley { \cite[Theorem 2.6]{St3}}] Let $G$ be a graph with vertex set $ [n]$. Then  
\begin{equation} \label{stan}
\omega X_G(\x)=\sum_{\pi \in \Pi_{G,n}}|\mu_G(\hat{0},\pi)|p_{\pa(\pi)}.
\end{equation}
Consequently $\omega X_G(\x)$ is $p$-positive.
\end{thm}

A poset $P$ is called {\em $(r+s)$-free} if it contains no induced subposet isomorphic to the direct sum  of an $r$ element chain and an $s$ element chain. 
A classical result (see \cite{ScSu}) says that a poset  is  a unit interval order if and only if  it is both $(3+1)$-free and $(2+2)$-free.

\begin{con} [Stanley and Stembridge { \cite[Conjecture 5.5]{StSte}, \cite[Conjecture 5.1]{St3}}]   \label{stancon} Let $G$ be the incomparability graph of a $(3+1)$-free poset.  Then $X_G(\x)$ is $e$-positive.
\end{con}

This conjecture is still open even for unit interval orders.  Gasharov \cite{Ga} obtains the weaker, but still very interesting, result that, under the assumptions of Conjecture \ref{stancon}, $X_G(\x)$ is Schur-positive. 
Let $P$ be a poset on $n$ and let $\lambda$ be a partition if $n$.   Gasharov defines a {\em $P$-tableau  of shape $\lambda$} to be a filling  of a Young diagram of shape $\lambda$ (in English notation)
with elements of $P$ such that
\begin{itemize}
  \item each element of $P$ appears exactly once,
  \item if $y \in P$ appears immediately to the right of $x \in P$, then $y>_P x$, and
  \item if $y \in P$ appears immediately below $x \in P$, then $y \not <_P x$.
\end{itemize}

Given a $P$-tableau $T$, let $\lambda(T)$ be the shape of $T$.  Let ${\mathcal T}_P$ be the set of all $P$-tableaux.

\begin{thm}[Gasharov \cite{Ga}] \label{gash}
Let $P$ be a $(3+1)$-free poset.  Then
\[
X_{\inc(P)}(\x)=\sum_{T \in {\mathcal T}_P} s_{\lambda(T)}.
\]
\end{thm}

\subsection{A quasisymmetric refinement}
We define the {\em chromatic quasisymmetric function} of $G$ as
$$X_G(\x,t) := \sum _{c\in C(G)} t^{\asc(c)} \x_c,$$
where 
$$\asc(c) :=| \{\{i,j\} \in E(G) : i < j \mbox{ and  } c(i) < c(j) \}|.$$
It is straightforward to confirm that $X_G(\x,t) \in \cq_{\zz}[t]$.

While the chromatic quasisymmetric function is defined for an arbitrary graph on vertex set $[n]$, our results will concern incomparability graphs of partially ordered sets.  The natural unit interval orders are particularly significant in our theory of chromatic quasisymmetric functions because they yield symmetric functions.
\begin{prop}[Shareshian and Wachs \cite{ShWa3}]  \label{symprop} Let $P$ be a natural unit interval order.   Then $$X_{\inc(P)}(\x,t) \in \Lambda_{\zz}[t].$$
\end{prop}

Before looking at some examples of $X_{\inc(P)}(\x,t)$ for natural unit interval orders, let us consider the poset 
$P$ on vertex set $[3]$ whose only relation is $1<2$.  Clearly $P$ is not a natural unit interval order.   (Recall $A_P(q,t)$ is given in (\ref{nonnatint}).)      Using Theorem~\ref{main} below we compute $$X_{\inc(P)}(\x,t)= (e_3 + F_{3,\{2\}}) + 2 e_3 t + (e_3 + F_{3,\{1\}})t^2,$$
which is clearly not in $\Lambda_\Z[t]$.

\begin{eg} $P=P_{n,1}$.  The incomparability graph has no edges. Hence 
\bq \label{pn1} X_{\inc(P_{n,1})}(\x,t)= e_1^n . \eq
\end{eg}
\begin{eg} $P=P_{n,n}$. The incomparability graph is the complete graph.  Each proper coloring $c$ is an injective map which can be associated with a permutation $\sigma \in \S_n$ for which $\inv(\sigma) = \asc(c)$.  It follows that
\bq \label{pn2} X_{\inc(P_{n,n})}(\x,t)= e_n  \sum_{\sigma \in \S_n} t^{\inv(\sigma)} = [n]_t! \,e_n.\eq
\end{eg}

\begin{eg} $P=P_{n,2}$.
 Recall that the incomparability graph  of  $P_{n,2}$ is a path.    To each proper coloring $c$ of  the path $\inc(P_{n,2})$ one  can associate the word $w(c):=c(n),c(n-1),\dots,c(1)$.  This word is clearly a Smirnov word of length $n$ (c.f. Section~\ref{smir}) and  $\des(w(c)) = \asc(c)$.  Since $w$ is a bijection from $C(\inc(P_{n,2}))$  to $W_n$, we have $$X_{\inc(P_{n,2})}(\x,t) = \sum_{j=0}^{n-1} Y_{n,j}(\x) t^j .$$
It therefore follows from Corollary~\ref{stancor} that 
\bq \label{torchrom} \omega X_{\inc(P_{n,2})}(\x,t) =  \sum_{j=0}^{n-1}  \ch H^{2j}(\v(\Delta_n)) t^j.\eq
\end{eg}

The following result is a consequence of Proposition~\ref{symprop}.

\begin{thm}[Shareshian and Wachs \cite{ShWa3}] \label{palinchromth} Let $P$ be a natural unit interval order.   Then  $X_{\inc(P)}(\x,t)$  is a palindromic polynomial in $t$.
\end{thm}

From (\ref{pn1}), (\ref{pn2}), (\ref{torchrom}) and Theorem~\ref{hunimodal}  we see that the following conjectured
 refinement of the unit interval order case of the Stanley-Stembridge conjecture (Conjecture~\ref{stancon}) is true for $P=P_{n,k}$ when $k=1,2,n$.

 \begin{con} \label{quasistan} Let $P$ be a natural unit interval order.   Then  $X_{\inc(P)}(\x,t)$  is an $e$-positive and $e$-unimodal polynomial in $t$. \end{con}
 
 We have verified $e$-positivity and $e$-unimodality for several other cases including
$P=P_{n,k}$ when  $k=n-1, n-2$, and by computer for all $n \le 8$.
Our computation yields,
\begin{eqnarray} \label{comps}\\
\nonumber X_{\inc(P_{n,n-1})}(\x,t) &=& [n-2]_t ! \left([n]_t[n-2]_t e_n + t^{n-2} e_{n-1,1}\right)\\ \label{comps2} \\
\nonumber X_{\inc(P_{n,n-2})}(\x,t) &=& [n-4]_t! \,\, ([n]_t[ n-3]_t^3 e_n  + [n-2]_t t^{n-3} ([n-3]_t +\\ 
\nonumber & &    [2]_t [n-4]_t) e_{n-1,1} + t^{2n-7}[2]_t e_{n-2,2}),
\end{eqnarray}
from which the conjecture is easily verified using Propositions~\ref{tooluni} and~\ref{tooluni2}.

In \cite{St3} Stanley proves that for any graph $G$ on $[n]$, the number of acyclic orientations of $G$ 
with $j$ sinks is equal to  $\sum_{\lambda \in \Par(n,j)}  c^G_\lambda$, where $\Par(n,j)$ is the set of partitions of $n$ into $j$ parts and  $c^G_
\lambda$ is the coefficient of $e_\lambda$   in the $e$-basis expansion of the chromatic symmetric function 
$X_G(\x)$. We obtain the following refinement of this result using essentially the same proof as Stanley's, thereby providing a bit of further evidence for $e$-positivity of $X_G(\x,t)$.

\begin{thm}[Shareshian and Wachs \cite{ShWa3}] \label{acyclicth} Let $G$ be the incomparability graph of   a natural unit interval order on $[n]$.   For each $\lambda \vdash n$, let $c^G_{\lambda}(t)$ be the coefficient of $e_\lambda$ in  the $e$-basis expansion of
$X_{G}(\x,t)$.  Then 
$$\sum_{\lambda \in \Par(n,j)}  c^G_\lambda(t) = \sum_{o \in {\mathcal O}(G,j)} t^{\asc(o)},$$
where ${\mathcal O}(G,j)$ is the set of acyclic orientations of $G$ with $j$ sinks and $\asc(o)$ is   the number of directed edges $(i,j)$  of $o$ for which $i<j$.  
\end{thm}

Theorem~\ref{acyclicth} gives a combinatorial description of  the coefficient of $e_n$ in the $e$-basis expansion of $X_G(\x,t)$.  In Theorem~\ref{ecoefth} and Conjecture~\ref{ecoefcon} below we give  alternative descriptions.  
Let $P$ be a poset on $[n]$.  We say that $\sigma \in \S_n$ has   a  left-to-right $P$-maximum at $r\in [n]$ if $\sigma(r) >_P \sigma(s) $ for $1\le s <r$.  The  left-to-right $P$-maximum at $1$ will be referred to as   trivial.  (The notion of (trivial) left-to-right $P$-maximum can be extended in an obvious way to permutations of any subset of $P$.) Let 
\begin{equation}  \label{defceq} c_P(t) :=\sum_{\sigma} t^{\inv_{\inc (P)}(\sigma)},\end{equation}
 where $\sigma$ ranges over the permutations in $\S_n$ with no $P$-descents and no nontrivial  left-to-right $P$-maxima. 

\begin{thm}[Shareshian and Wachs \cite{ShWa3}] \label{ecoefth} Let $P$ be a natural unit interval order on $[n]$.  The coefficient of $e_n$ in the expansion of $X_{\inc(P)}(\x,t)$ in the $e$-basis of $\Lambda_{\Z[t]}$ is equal to $c_P(t)$. \end{thm}

\begin{con} \label{ecoefcon} Let $P$ be a natural unit interval order on $[n]$.  Then $$c_P(t) = [n]_t \prod_{i =2}^n [a_i]_t ,$$
where $a_i = |\{j \in [i-1]: \{i,j\} \in E(\inc(P))\}|$.   Consequently $c_P(t)$ is palindromic and unimodal. \end{con}

In the case that $P=P_{n,2}$, this conjecture is true by Theorem~\ref{ecoefth} and equations (\ref{formQ}) and (\ref{torchrom}).
It is also true for $P_{n,k}$ when $k=1,n, n-1,n-2$ by Theorem~\ref{ecoefth} and equations   (\ref{pn1}), (\ref{pn2}), (\ref{comps}) and (\ref{comps2}), respectively.

\subsection{Schur and power sum decompositions} When $P$ is a natural unit interval order, we have the following  refinement of Gasharov's Schur positivity result (Theorem~\ref{gash}).   For a $P$-tableau $T$ and a graph $G$ with vertex set $[n]$, let $\inv_G(T)$ be the number of edges $\{i,j\} \in E(G)$ such that $i<j$ and $i$ appears to the south of $j$ in $T$. 

\begin{thm}[Shareshian and Wachs \cite{ShWa3}]\label{schurcon}
Let $P$ be a natural unit interval order poset on $[n]$ and
 let $G$ be the incomparability graph of $P$.  
  Then
\[
X_G(\x,t)=\sum_{T \in {\mathcal T}_P}t^{\inv_G(T)}s_{\lambda(T)}.
\]
\end{thm}

Now by (\ref{torchrom}) we have the following.
\begin{cor} \label{swdecomp} Let $\lambda \vdash n$.  The multiplicity of the irreducible Specht module $S^\lambda$ in $H^{2j}(\v(\Delta_n)) $ is equal to the number of $P_{n,2}$-tableaux of shape $\lambda^\prime$ with $\inv_{\inc(P_{n,2})}(T) = j$.
\end{cor}

By comparing this decomposition to Stembridge's decomposition (Theorem~\ref{stemdecomp}) we see that the number of marked tableaux of shape $\lambda$ and index $j$ equals the number of $P_{n,2}$-tableaux of shape $\lambda^\prime$ and $\inv_{\inc(P_{n,2})}(T) = j$.  It would be interesting to find   a bijective proof of this fact.

 Theorem \ref{schurcon} shows that when $G$ is the incomparability graph of a natural unit interval order, the coefficient of each power of $t$ in $\omega X_G(\x,t)$ is a nonnegative integer combination of Schur functions and therefore the Frobenius characteristic of an actual representation of $\S_n$.  Conjecture~\ref{quasistan} says that this linear representation arises from a permutation representation in which each point stabilizer is a Young subgroup.  In Section~\ref{hesssec} we present a very promising concrete candidate for the desired permutation representation.
 
 Next we attempt to refine Stanley's  $p$-basis decomposition of the chromatic symmetric functions  (Theorem~\ref{stan}).  In Conjectures \ref{powercon1} and \ref{powercon2} below, we provide two proposed formulae for this power sum decomposition.
  With $\mu=(\mu_1 \ge \mu_2 \ge \dots \ge \mu_l)$ a partition of $n$ and $P$ a natural unit interval order on $[n]$, we call a permutation $\sigma \in \S_n$, {\em $(P,\mu,1)$-compatible} 
 if, when we break $\sigma$ (in one line notation) into consecutive segments of lengths $\mu_1,\ldots,\mu_l$, the segments have no $P$-descents and no nontrivial left-to-right $P$-maxima.  
We call $\sigma \in \S_n$, {\em $(P,\mu,2)$-compatible} if, 
when we break $\sigma$  into consecutive segments of lengths $\mu_1,\ldots,\mu_l$, the segments have no $P$-ascents and  they begin with the numerically smallest letter of the segment.
We write $\S_{P,\mu,i}$, where $i=1,2$, for the set of all $(P,\mu,i)$-compatible elements of $\S_n$.
 \begin{con}  \label{powercon1} Let $P$ be a natural unit interval order on $[n]$ with incomparability graph $G$.  Then
 $$\omega X_G(\x,t) = \sum_{\mu \vdash n}  z_\mu^{-1} p_\mu  \sum_{\sigma \in \S_{P,\mu,1}} t^{\inv_G(\sigma)}.$$
  \end{con}
 
  \begin{con}   \label{powercon2} Let $P$ be a natural unit interval order on $[n]$ with incomparability graph $G$.  Then
  $$\omega X_G(\x,t) = \sum_{\mu \vdash n}  z_\mu^{-1} p_\mu \,\,\prod_{i=1}^{l(\mu)} [\mu_i]_t\sum_{\sigma \in \S_{P,\mu,2}} t^{\inv_G(\sigma)}.$$
 \end{con}
 
It can be shown that both conjectures reduce to the formula in Theorem~\ref{stan} when $t=1$.
Conjecture~\ref{powercon1} says that the coefficient of $z_{(n)}^{-1} p_{(n)}$ is  equal to 
$c_P(t)$ defined in  (\ref{defceq}).   Since the coefficient of $h_{(n)}$ in the $h$-basis 
decomposition must equal the coefficient of $ z_{(n)}^{-1} p_{(n)}$ in the $p$-basis 
decomposition, by Theorem~\ref {ecoefth}, Conjecture~\ref{powercon1} gives the correct coefficient of $
 z_{(n)}^{-1} p_{(n)}$.
When $P=P_{n,2}$ it is not difficult to see that both conjectures reduce to Theorem~\ref{pdecomp}.

\subsection{ Fundamental quasisymmetric function basis decomposition} \label{fundquasec}
\begin{thm}[Shareshian and Wachs \cite{ShWa3}] \label{main} Let  $G$ be the incomparability graph of a  poset $P$ on $[n]$.   Then 
$$\omega X_G(\x,t) = \sum_{\sigma \in \S_n} t^{\inv_G(\sigma)} F_{n,\Des_{P}(\sigma)},$$
where $\omega$ is the involution on $\cq_{\Z}$ that maps $F_{n,S}$ to $F_{n,[n-1] \setminus S}$ for each $n \in \nn$ and $S \subseteq [n-1]$.  (This extends  the involution $\omega$ on $\Lambda_{\Z}$ that maps $h_n$ to $e_n$)
\end{thm}

Theorem \ref{main} refines Corollary 2 in Chow's paper \cite{Ch}.  Indeed, one obtains Chow's result by setting $t=1$ in Theorem \ref{main}.  Our proof of Theorem \ref{main} follows the same path as Chow's proof of his Corollary~2.

We use Theorem~\ref{main} and (\ref{sps}) to compute the principal stable specialization of $\omega X_G(\x,t) $.
\begin{cor} \label{maincor} Let  $G$ be the incomparability graph of a  poset $P$ on $[n]$. Then
$$(q;q)_n \Ps( \omega X_G(\x,t)) = A_P(q,t).$$
\end{cor}

Recall that for $P=P_{n,2}$ we have a  formulation for the  expansion of $\omega X_{\inc(P)}(\x,t)$ in the fundamental quasisymmetric function basis that is different from  that of Theorem~\ref{main}.  It is obtained by combining (\ref{eulfund}) and (\ref{carlg}).  By equating these formulations, we obtain the  identity,
\begin{equation}\label{fundeq}\sum_{\sigma \in \S_n} t^{\des(\sigma^{-1})} F_{n,\Des_{\ge 2}(\sigma)} = \sum_{\sigma \in \S_n} t^{\exc(\sigma)} F_{n,\Dex(\sigma)} .\end{equation}
Taking the stable principle specialization of both sides of (\ref{fundeq}) yields the following new Euler-Mahonian result.

\begin{thm}[Shareshian and Wachs \cite{ShWa3}] \label{eqEu} For all $n\in \pp$,
$$ \sum_{\sigma \in \S_n} q^{\rmaj_{2} (\sigma)} t^{ \des (\sigma^{-1}) } = \sum_{\sigma \in \S_n} q^{\maj(\sigma)}t^{\exc(\sigma)} .$$
\end{thm}
It would be interesting to find bijective proofs of ({\ref{fundeq}) and Theorem~\ref{eqEu}.

Theorem~\ref{palinAth}, asserting  palindromicity of $A_P(q,t)$, is  proved by taking  stable principal specializations in  Theorem~\ref{palinchromth}, by means of Corollary~\ref{maincor}.  Conjecture~\ref{unieucon}, asserting  unimodality of $A_P(q,t)$,  can be obtained by taking  stable principal specializations in  Conjecture~\ref{quasistan} (see Lemma~\ref{hposqpos}).  By taking a specialization called (nonstable) principal specialization, a stronger unimodality conjecture follows from Conjecture~\ref{quasistan}.  For a poset $P$ on $[n]$ and $\sigma \in \S_n$, define $\des_P(\sigma) := |\Des_P(\sigma)|$.

\begin{con} \label{triple} Let $P$ a natural unit interval order on $[n]$. Then 
$$A_P(q,p,t):=\sum_{\sigma \in \S_n} q^{\maj_P(\sigma)} p^{\des_P(\sigma)} t^{\inv_{\inc(P)}(\sigma)}$$ is a palindromic unimodal polynomial in $t$.
\end{con}

It follows from Corollary~\ref{maincor} and Lemma~\ref{thdes}  that 
$A_P(\xi_d,t)$ is the coefficient of $z_{d^m}^{-1}p_{d^m}$ in the $p$-basis expansion of $\omega X_{\inc(P)}(\x,t)$ when $P$ is a natural unit interval order on $[dm]$. Hence Conjecture~\ref{genunity} on the evaluation of $A_P(q,t)$ at roots of unity is a consequence of Conjectures~\ref{quasistan} and~\ref{powercon2}.

\section{Hessenberg varieties} \label{hesssec}

 Let $G=GL_n(\cc)$ and let $B$ be the set of upper triangular matrices in $G$.  The (type A) {\em flag variety} is the quotient space $G/B$. 
 Fix now a natural unit interval order $P$ and let $M_{n,P}$ be the set of all $n \times n$ complex matrices $(a_{ij})$ such that $a_{ij}=0$ whenever $i >_P j$.   Fix a nonsingular diagonal $n \times n$ matrix $\bft$ with $n$ distinct eigenvalues.  The {\em regular semisimple Hessenberg variety of type $A$} associated to $P$ is
\[
\hess (P):=\{gB \in G/B:g^{-1}\bft g \in M_{n,P}\}.
\]
(Note that $\hess (P)$ is well defined, since the group $B$ normalizes the set $M_{n,P}$.)  

Certain regular semisimple Hessenberg varieties of type A were studied initially by De Mari and Shayman in \cite{DeSh}.  Hessenberg varieties for other Lie types are defined and studied by De Mari, Procesi and Shayman in \cite{DePrSh}.  Such varieties are determined by certain subsets of a root system.  In \cite[Theorem 11]{DePrSh} it is noted that, for arbitrary Lie type, the Hessenberg variety associated with the set of simple roots is precisely the  toric variety associated with the corresponding Coxeter complex.  
In particular, in Lie type A, the poset giving rise to the regular semisimple Hessenberg variety associated with simple roots is $P_{n,2}$ since $E(\inc(P_{n,2}))=\{\{i,i+1\}: i \in [n-1]\}$. Thus
\begin{equation}\label{torhess}  \hess(P_{n,2}) = \v(\Delta_n).\end{equation}

De Mari, Procesi and Shayman also show that the cohomology of $\hess (P)$ is concentrated in even degrees  and that for $0\le j \le   |E(\inc(P))|$, 
\begin{equation} \label{dmpseq} \dim H^{2j}(\hess(P)) = a_{P,j}(1).\end{equation}
Hence unimodality of the $P$-Eulerian polynomials $A_P(1,t)$ (Theorem \ref{hesseuth}) follows from  the hard Lefschetz theorem since $\hess(P)$ is smooth.
Our approach to establishing the unimodality of the $(q,P)$-Eulerian polynomial $A_P(q,t)$ is to find a representation of the symmetric group on $H^{2j}(\hess(P))$ such that the stable principal specialization  of its Frobenius characteristic yields $a_{P,j}(q)$.  If the hard Lefschetz map commutes with the action of the symmetric group then it follows from Schur's lemma and Lemma~\ref{hposqpos} that $A_P(q,t)$ is unimodal.

For the cases $P= P_{n,k}$, where k=1,2, we already have the desired representations.   Indeed,  $\hess(P_{n,1})$ consists of $n!$ isolated points and the representation of $\S_n$ on $H^0(\hess(P_{n,1}))$ is the regular representation.  
For $k=2$, we can use the representation of $\S_n$ on $H^*(\hess(P_{n,2}))= H^*(\v(\Delta_n))$  discussed in Section~\ref{eultorsec}.

 For general natural unit interval orders $P$, Tymoczko   \cite{Ty3}, \cite{Ty2} has defined a representation of $\S_n$ on $H^*(\hess(P))$ via a theory of Goresky, Kottwitz and MacPherson known as GKM theory (see \cite{GoKoMacP}, and see \cite{Ty} for an introductory description of GKM theory).
The centralizer $T=C_G(\bft)$ consists of the diagonal matrices in $G:=GL_n(\cc)$.  It follows that the torus $T$ acts (by left translation) on $\hess (P)$.  The technical conditions required for application of GKM theory  are satisfied by this action, and it follows that one can describe the cohomology of $\hess(P)$ using the moment graph associated to this action.  The moment graph $M$ is a subgraph of the Cayley graph of $\S_n$ with generating set consisting of  the transpositions.  The vertex set of $M$ is $\S_n$ and the edges connect pairs of elements that differ by a transposition $(i,j)$ such that $\{i,j\} \in E(\inc(P)) $.
Thus $M$ admits an action of $\S_n$, and this action can be used to define a linear representation of $\S_n$ on the cohomology $H^\ast(\hess (P))$.  

 Tymoczko's representation of $\S_n$ on $H^*(\hess(P))$ in  the cases $P =P_{n,k}$,  $k=1,2$,  is the same as the respective representations of  $\S_n$  discussed above.  In the case that $P=P_{n,n}$, it follows from \cite[Proposition 4.2]{Ty2}  that  for all $j$, Tymoczko's representation  is  isomorphic to $a^{(n)}_{n,j}$ copies of the trivial representation.
 
By applying the hard Lefschetz theorem,  MacPherson and Tymoczko obtain the following result.

\begin{thm} [MacPherson and Tymoczko \cite{MacTy}] \label{MTschur}
For all natural unit interval orders $P$, the sequence $(\ch H^{2j}(\hess(P)))_{0 \le j \le |E(\inc(P)|}$ is palindromic and Schur-unimodal.
\end{thm}

Tymoczko poses the following problem.

\begin{prob}[Tymoczko \cite{Ty2}] \label{tymprob}  Given any natural unit interval order $P$ on vertex set $[n]$, describe the decomposition of the representation of $\S_n$ on $H^{2j}(\hess(P))$ into irreducibles.
\end{prob}

We finally come to the conjecture that ties together the three topics of this paper.

\begin{con}\label{hesschrom} Let $P$ be a natural unit interval order on $[n]$.
Then \begin{equation}\label{mainconeq}\sum_{j=0}^{|E(\inc(P))|} \ch H^{2j}(\hess(P)) t^j = \omega X_{\inc(P)}(\x,t).\end{equation}
\end{con}

By (\ref{torhess}) and (\ref{torchrom}) the conjecture is true in the case that $P=P_{n,2}$.
It is straightforward to verify in the case $P=P_{n,k}$,  when $k=1,n-1,n$, and for all natural unit interval orders on $[n]$ when $n \le 4$.  The conjecture is also true for the  parabolic Hessenberg varieties studied by Teff \cite{Te}.    It follows from (\ref{dmpseq}) that the coefficient of the monomial symmetric function $m_{1^n}$ in the expansion of the left hand side of (\ref{mainconeq}}) in the monomial symmetric function basis equals the coefficient of $m_{1^n}$ in the expansion of the right hand side of (\ref{mainconeq}).  

Conjecture~\ref{hesschrom}, if true,  would have many important ramifications.  It would allow us to transfer what we know about chromatic quasisymmetric functions to Tymoczko's representation and vice-versa.   For instance Theorem~\ref{schurcon} would provide a solution to Tymoczko's problem. Conjecture~\ref{hesschrom}   provides an approach to proving the Stanley-Stembridge conjecture (Conjecture~\ref{stancon}) for unit interval orders.   Indeed  one would only have to prove the following conjecture. 

\begin{con} \label{tympermcon} For all natural interval orders $P$, Tymoczko's representation of $\S_n$ on $H^*(\hess(P))$ is a permutation representation  in which each point stabilizer is a Young subgroup.
\end{con} 

   Our refinement of the Stanley-Stembridge conjecture (Conjecture~\ref{quasistan}) would be equivalent to the following strengthening of the result of MacPherson and Tymoczko.

\begin{con}  \label{SWhess}
For all natural unit interval orders $P$, the palindromic sequence $(\ch H^{2j}(\hess(P)))_{0 \le j \le |E(\inc(P)|}$ is  h-unimodal.
\end{con}

To establish the  unimodality of the $(q,P)$-Eulerian numbers (Conjecture~\ref{unieucon}) one would only need to prove Conjecture~\ref{hesschrom}.  Indeed  Conjecture~\ref{hesschrom}, Theorem~\ref{MTschur}, Corollary~\ref{maincor} and Lemma~\ref{hposqpos} together imply the unimodality result.   The more general unimodality conjecture for $A_P(q,p,t)$ (Conjecture~\ref{triple}) is also a consequence of Conjecture~\ref{hesschrom}.   

\section*{Acknowledgments}  We are grateful to Richard Stanley for informing us of Theorem~\ref{stanth}, which led us to  introduce the  chromatic quasisymmetric functions.   
We also  thank Julianna Tymoczko for   useful
discussions and the anonymous referee for helpful comments.

\end{document}